# WEIGHTED EMPIRICAL LIKELIHOOD IN SOME TWO-SAMPLE SEMIPARAMETRIC MODELS WITH VARIOUS TYPES OF CENSORED DATA

By Jian-Jian Ren[1]

*University of Central Florida*

In this article, the *weighted empirical likelihood* is applied to a general setting of two-sample semiparametric models, which includes biased sampling models and case-control logistic regression models as special cases. For various types of censored data, such as right censored data, doubly censored data, interval censored data and partly interval-censored data, the weighted empirical likelihood-based semiparametric maximum likelihood estimator $(\tilde{\theta}_n, \tilde{F}_n)$ for the underlying parameter $\theta_0$ and distribution $F_0$ is derived, and the strong consistency of $(\tilde{\theta}_n, \tilde{F}_n)$ and the asymptotic normality of $\tilde{\theta}_n$ are established. Under biased sampling models, the weighted empirical log-likelihood ratio is shown to have an asymptotic *scaled* chi-squared distribution for censored data aforementioned. For right censored data, doubly censored data and partly interval-censored data, it is shown that $\sqrt{n}(\tilde{F}_n - F_0)$ weakly converges to a centered Gaussian process, which leads to a consistent goodness-of-fit test for the case-control logistic regression models.

**1. Introduction.** Consider the following two-sample semiparametric model:

(1.1)
$X_1, \ldots, X_{n_0}$ is a random sample with density $f_0(x)$,

$Y_1, \ldots, Y_{n_1}$ is a random sample with density $g_0(x) = \varphi(x; \theta_0) f_0(x)$,

where the two samples are independent, and $\varphi(x; \theta_0)$ is a known function with $x \in \mathbb{R}$ and a unique unknown parameter $\theta_0 \in \mathbb{R}^q$, while $f_0$ and $g_0$ are the density functions of unknown nonnegative distribution functions (d.f.)

---

Received September 2005; revised February 2007.
[1]Research supported in part by NSF Grants DMS-02-04182 and DMS-06-04488.
*AMS 2000 subject classifications.* 62N02, 62N03, 62N01.
*Key words and phrases.* Biased sampling, bootstrap, case-control data, doubly censored data, empirical likelihood, Kolmogorov–Smirnov statistic, interval censored data, likelihood ratio, logistic regression, maximum likelihood estimator, partly interval-censored data, right censored data.

$F_0$ and $G_0$, respectively. This model (1.1) includes biased sampling models (Vardi [32]) and case-control logistic regression models (Prentice and Pyke [22]) as special cases, for which there has not been any published work dealing with censored data. In this article, we study model (1.1) when at least one of the two samples is not completely observable due to censoring. As follows, we use random sample $X_1, \ldots, X_{n_0}$ to illustrate the censoring models under consideration here, while Examples 1 and 2 discuss biased sampling models and case-control logistic regression models, respectively.

*Right censored sample.* The observed data are $\boldsymbol{O}_i = (V_i, \delta_i), 1 \leq i \leq n_0$, with

(1.2) $\quad V_i = \begin{cases} X_i, & \text{if } X_i \leq C_i, \quad \delta_i = 1, \\ C_i, & \text{if } X_i > C_i, \quad \delta_i = 0, \end{cases}$

where $C_i$ is the right censoring variable and is independent of $X_i$. This type of censoring has been extensively studied in the literature in the past few decades.

*Doubly censored sample.* The observed data are $\boldsymbol{O}_i = (V_i, \delta_i)$, $1 \leq i \leq n_0$, with

(1.3) $\quad V_i = \begin{cases} X_i, & \text{if } D_i < X_i \leq C_i, \quad \delta_i = 1, \\ C_i, & \text{if } X_i > C_i, \quad \delta_i = 2, \\ D_i, & \text{if } X_i \leq D_i, \quad \delta_i = 3, \end{cases}$

where $C_i$ and $D_i$ are right and left censoring variables, respectively, and they are independent of $X_i$ with $P\{D_i < C_i\} = 1$. This type of censoring has been considered by Turnbull [31], Chang and Yang [4], Gu and Zhang [11] and Mykland and Ren [17], among others. One recent example of doubly censored data was encountered in a study of primary breast cancer (Ren and Peer [28]).

*Interval censored sample.*

CASE 1.  The observed data are $\boldsymbol{O}_i = (C_i, \delta_i)$, $1 \leq i \leq n_0$, with

(1.4) $\quad\quad\quad\quad\quad\quad \delta_i = I\{X_i \leq C_i\}.$

CASE 2.  The observed data are $\boldsymbol{O}_i = (C_i, D_i, \delta_i)$, $1 \leq i \leq n_0$, with

(1.5) $\quad\quad\quad\quad \delta_i = \begin{cases} 1, & \text{if } D_i < X_i \leq C_i, \\ 2, & \text{if } X_i > C_i, \\ 3, & \text{if } X_i \leq D_i, \end{cases}$

where $C_i$ and $D_i$ are independent of $X_i$ and satisfy $P\{D_i < C_i\} = 1$ for Case 2. These two types of interval censoring were considered by Groeneboom and Wellner [10], among others. In practice, interval censored Case 2 data were encountered in AIDS research (Kim, De Gruttola and Lagakos [16]; see discussion in Ren [26]).

*Partly interval-censored sample.*

"CASE 1" PARTLY INTERVAL-CENSORED DATA. The observed data are

$$(1.6) \qquad \boldsymbol{O}_i = \begin{cases} X_i, & \text{if } 1 \le i \le k_0, \\ (C_i, \delta_i), & \text{if } k_0 + 1 \le i \le n_0, \end{cases}$$

where $\delta_i = I\{X_i \le C_i\}$ and $C_i$ is independent of $X_i$.

GENERAL PARTLY INTERVAL-CENSORED DATA. The observed data are

$$(1.7) \qquad \boldsymbol{O}_i = \begin{cases} X_i, & \text{if } 1 \le i \le k_0, \\ (\boldsymbol{C}, \boldsymbol{\delta}_i), & \text{if } k_0 + 1 \le i \le n_0, \end{cases}$$

where for $N$ potential examination times $C_1 < \cdots < C_N$, letting $C_0 = 0$ and $C_{N+1} = \infty$, we have $\boldsymbol{C} = (C_1, \ldots, C_N)$ and $\boldsymbol{\delta}_i = (\delta_i^{(1)}, \ldots, \delta_i^{(N+1)})$ with $\delta_i^{(j)} = 1$, if $C_{j-1} < X_i \le C_j$; 0, elsewhere. This means that for intervals $(0, C_1], (C_1, C_2], \ldots, (C_N, \infty)$, we know in which one of them $X_i$ falls. These two types of partly interval-censoring were considered by Huang [12], among others. As pointed out by Huang [12], in practice the general partly interval-censored data were encountered in Framingham Heart Disease Study (Odell, Anderson and D'Agostino [18]), and in the study on incidence of proteinuria in insulin-dependent diabetic patients (Enevoldsen et al. [5]).

EXAMPLE 1 (*Biased sampling model*). In (1.1), let

$$(1.8) \qquad \varphi(x; \theta_0) = \theta_0 w(x), \qquad \theta_0 \in \mathbb{R},$$

where $w(x)$ is a *weight function* with positive value on the support of $F_0$, and $\theta_0 = 1/w_0$ is the *weight parameter* satisfying $w_0 = \int_0^\infty w(x) \, dF_0(x)$. Then, (1.1) is a *two-sample biased sampling problem*, for which the case with *length-biased* distribution $G_0$, that is, $w(x) = x$ in (1.8), was considered by Vardi [32], and the empirical log-likelihood ratio for the mean of $F_0$ was shown to have an asymptotic chi-squared distribution by Qin [23]. More general biased sampling models were considered by Vardi [33], Gill, Vardi and Wellner [9], who discussed various application examples, and showed that the maximum likelihood estimator for $F_0$ is asymptotically Gaussian and efficient. For right censored samples in (1.1), Vardi [33] gave an estimator for $F_0$ based on the EM algorithm, but the asymptotic properties of the estimator were not studied. Below, we discuss practical examples of biased sampling problem with censored data.

In Patil and Rao [20], the biased sampling problem is discussed in the context of efficiency of early screening for disease. Using our notations in (1.1), if $F_0$ is the d.f. of the duration of the preclinical state of certain chronic disease, then the first sample in (1.1) is taken from those whose clinical state is detected by the usual medical care. If at a certain point in time some individuals in the preclinical state begin participating in an early detection program, then such a program identifies them by a length-biased sampling. In other words, the second sample in (1.1) is taken from those who participated in the early detection program, and $G_0$ is a length-biased distribution. However, in reality a usual screening program for "disease" is conducted by examining an individual periodically with a fixed length of time between two consecutive check-ups. The data encountered in such a screening program is typically a doubly censored sample (1.3); that is, the actually observed data for the second sample in (1.1) is doubly censored. In statistical literature, examples of doubly censored data encountered in screening programs have been given by Turnbull [31] and Ren and Gu [27], among others.

EXAMPLE 2 (*Case-control logistic regression model*). In (1.1), let

$$\varphi(x; \theta_0) = e^{\alpha_0 + \beta_0 x},$$
$$F_0(x) = P\{T \le x | Z = 0\}, \qquad G_0(x) = P\{T \le x | Z = 1\}; \tag{1.9}$$

then under reparameterization by Qin and Zhang [24], model (1.1) is equivalent to the following *case-control logistic regression model* (Prentice and Pyke [22]):

$$P\{Z = 1 | T = x\} = \frac{\exp(\alpha^* + \beta_0 x)}{1 + \exp(\alpha^* + \beta_0 x)}, \tag{1.10}$$

where $\theta_0 = (\alpha_0, \beta_0) \in \mathbb{R}^2$, $Z$ is the binary response variable (with value 1 or 0 to indicate presence or absence of a disease or occurrence of an event of interest), $T$ is the covariate variable, and $(\alpha^*, \beta_0)$ is the regression parameter satisfying $\alpha_0 = \alpha^* + \ln[(1 - \pi)/\pi]$ for $\pi = P\{Z = 1\}$. Qin and Zhang [24] established asymptotic normality of the *semiparametric maximum likelihood estimators* (SPMLE) for $(\theta_0, F_0)$ in (1.9) with two complete samples in (1.1), and provided a goodness-of-fit test for (1.10). Below, we discuss an example to illustrate the situation with censored covariate variable $T$.

In the example of early detection of breast cancer considered by Ren and Gu [27], $T$ is the age at which the tumor could be detected when screening mammogram is the only detection method, and based on series screening mammograms the observed data on $T$ are doubly censored. This example is part of a study on the effectiveness of screening mammograms; see Ren and Peer [28] for precise description of left and right censored observations. Here, to study the effects of screening mammograms on survival, we consider those individuals who had breast cancer, and let $Z = 1$ represent death due to breast cancer within 5 years of diagnosis; $Z = 0$, otherwise. Then under (1.9), for those "dead" (i.e., $Z = 1$) the second sample in (1.1) is taken from the available screening mammogram records; thus the actually observed data from $G_0(x) = P\{T \leq x | Z = 1\}$ is a doubly censored sample. Similarly, for those "survived" (i.e., $Z = 0$) the first sample in (1.1) is also taken from screening mammogram records; thus also a doubly censored sample. Fitting the logistic regression model (1.10) with these two doubly censored case-control samples, we obtain $P\{Z = 1 | T = x_0\}$, which is the probability of "death" for an individual whose tumor was detected by screening mammogram at age $x_0$.

In this article, we apply *weighted empirical likelihood* (Ren [25]) to model (1.1) with the following two independent samples for $n = n_0 + n_1$:

(1.11)
$\boldsymbol{O}_1^X, \ldots, \boldsymbol{O}_{n_0}^X$ is the observed sample for sample $X_1, \ldots, X_{n_0}$,

$\boldsymbol{O}_1^Y, \ldots, \boldsymbol{O}_{n_1}^Y$ is the observed sample for sample $Y_1, \ldots, Y_{n_1}$,

where $\boldsymbol{O}_i^X$'s or $\boldsymbol{O}_j^Y$'s is possibly one of those censored samples described above, and we denote $\hat{F}$ and $\hat{G}$ as the *nonparametric maximum likelihood estimators* (NPMLE) for $F_0$ and $G_0$ based on $\boldsymbol{O}_i^X$'s and $\boldsymbol{O}_j^Y$'s, respectively. Section 2 provides a heuristic explanation of the concept of weighted empirical likelihood. For censored data (1.2)–(1.7) aforementioned, Section 3 derives the weighted empirical likelihood-based SPMLE $(\tilde{\theta}_n, \tilde{F}_n)$ for $(\theta_0, F_0)$, and establishes the strong consistency of $(\tilde{\theta}_n, \tilde{F}_n)$ and the asymptotic normality of $\tilde{\theta}_n$, while Section 4 further discusses Example 1 on biased sampling models, and shows that the weighted empirical log-likelihood ratio has an asymptotic scaled chi-squared distribution. For right censored data, doubly censored data and partly interval-censored data, Section 3 also shows that $\sqrt{n}(\tilde{F}_n - F_0)$ weakly converges to a centered Gaussian process, while Section 5 further discusses Example 2 on case-control logistic regression models, and provides a consistent goodness-of-fit test.

We note that the weighted empirical likelihood approach used in this article can be adapted to deal with more general biased sampling models. Also note that based on Ren and Gu [27], our results here on the case-control logistic regression models can be extended to $k$-dimensional ($k > 1$) covariate $\boldsymbol{T}$, where $\boldsymbol{T}$ contains one component that is subject to right censoring or doubly censoring.

For interval censored data (1.4)–(1.5), the weighted empirical likelihood approach enables us to obtain the strong consistency of the SPMLE $(\tilde{\theta}_n, \tilde{F}_n)$, the asymptotic normality of $\tilde{\theta}_n$, and the limiting distribution of the log-likelihood ratio via the asymptotic results on the NPMLE $\hat{F}$ or $\hat{G}$ for interval censored data by Groeneboom and Wellner [10] and Geskus and Groeneboom [6], among others. However, the techniques used in our proofs show that the weak convergence of $\tilde{F}_n$ for interval censored data relies on that of $\hat{F}$ or $\hat{G}$ for interval censored data, which is now unknown.

**2. Weighted empirical likelihood.** For random sample $X_1, \ldots, X_{n_0}$ from d.f. $F_0$, the empirical likelihood function (Owen [19]) is given by $L(F) = \prod_{i=1}^{n_0} [F(X_i) - F(X_i-)]$, where $F$ is any d.f. The weighted empirical likelihood function in Ren [25] may be understood as follows.

For each type of censored data aforementioned, the likelihood function has been given in literature, and the NPMLE $\hat{F}$ for $F_0$ is the solution which maximizes the likelihood function. Moreover, it is shown that from observed censored data $\{\boldsymbol{O}_i^X; 1 \leq i \leq n_0\}$, there exist $m_0$ distinct points $W_1^X < W_2^X < \cdots < W_{m_0}^X$ along with $\hat{p}_j^X > 0$, $1 \leq j \leq m_0$, such that $\hat{F}$ can be expressed as $\hat{F}(x) = \sum_{i=1}^{m_0} \hat{p}_i^X I\{W_i^X \leq x\}$ for above right censored data (Kaplan and Meier [15]), doubly censored data (Mykland and Ren [17]), interval censored data Case 1 and Case 2 (Groeneboom and Wellner [10]) and partly interval-censored data (Huang [12]). Since in all these cases $\hat{F}$ is shown to be a strong uniform consistent estimator for $F_0$ under some suitable conditions, we may expect a random sample $X_1^*, \ldots, X_{n_0}^*$ taken from $\hat{F}$ to behave asymptotically the same as $X_1, \ldots, X_{n_0}$. If $F_{n_0}^*$ denotes the empirical d.f. of $X_1^*, \ldots, X_{n_0}^*$, then from $\hat{F} \approx F_{n_0}^*$ we have

$$\prod_{i=1}^{n_0} P\{X_i = x_i\} \approx \prod_{i=1}^{n_0} P\{X_i^* = x_i^*\} = \prod_{j=1}^{m_0} (P\{X_1^* = W_j^X\})^{k_j}$$

$$\approx \prod_{j=1}^{m_0} (P\{X_1^* = W_j^X\})^{n_0[\hat{F}(W_j^X) - \hat{F}(W_j^X -)]}$$

$$= \prod_{j=1}^{m_0} (P\{X_1^* = W_j^X\})^{n_0 \hat{p}_j^X},$$

where $k_j = n_0[F_{n_0}^*(W_j^X) - F_{n_0}^*(W_j^X -)]$. Thus, the *weighted empirical likelihood function* (Ren [25])

$$(2.1) \qquad \hat{L}(F) = \prod_{i=1}^{m_0} [F(W_i^X) - F(W_i^X -)]^{n_0 \hat{p}_i^X}$$

may be viewed as the asymptotic version of the empirical likelihood function $L(F)$ for censored data. When there is no censoring, $\hat{L}(F)$ coincides with $L(F)$.

**3. SPMLE and asymptotic results.** This section derives the *semiparametric maximum likelihood estimator* for $(\theta_0, F_0)$ in (1.1) using censored data (1.11), and studies related asymptotic properties.

As general notations throughout this paper, let $\hat{F}$ and $\hat{G}$ be the NPMLE for $F_0$ and $G_0$ in (1.1) based on observed censored data $\boldsymbol{O}_1^X, \ldots, \boldsymbol{O}_{n_0}^X$ and $\boldsymbol{O}_1^Y, \ldots, \boldsymbol{O}_{n_1}^Y$ in (1.11), respectively. From Section 2, we know that there exist distinct points $W_1^X < \cdots < W_{m_0}^X$ and $W_1^Y < \cdots < W_{m_1}^Y$ with $\hat{p}_i^X > 0$ and $\hat{p}_i^Y > 0$ such that $\hat{F}$ and $\hat{G}$ can be expressed as

$$(3.1) \quad \hat{F}(x) = \sum_{i=1}^{m_0} \hat{p}_i^X I\{W_i^X \leq x\} \quad \text{and} \quad \hat{G}(x) = \sum_{i=1}^{m_1} \hat{p}_i^Y I\{W_i^Y \leq x\}$$

respectively, for those censored data aforementioned. We also let

$$(3.2) \quad \begin{aligned} (W_1, \ldots, W_m) &= (W_1^X, \ldots, W_{m_0}^X, W_1^Y, \ldots, W_{m_1}^Y), \\ (\hat{p}_1, \ldots, \hat{p}_m) &= (\hat{p}_1^X, \ldots, \hat{p}_{m_0}^X, \hat{p}_1^Y, \ldots, \hat{p}_{m_1}^Y), \\ (\omega_1, \ldots, \omega_m) &= (\rho_0 \hat{p}_1^X, \ldots, \rho_0 \hat{p}_{m_0}^X, \rho_1 \hat{p}_1^Y, \ldots, \rho_1 \hat{p}_{m_1}^Y), \end{aligned}$$

where $m = m_0 + m_1$, $\rho_0 = n_0/n$ and $\rho_1 = n_1/n$.

To derive an estimator for $(\theta_0, F_0)$ using both samples in (1.11), we apply weighted empirical likelihood function (2.1) to model (1.1), and obtain

$$\left( \prod_{i=1}^{m_0} [F(W_i^X) - F(W_i^X-)]^{n_0 \hat{p}_i^X} \right) \left( \prod_{j=1}^{m_1} [G(W_j^Y) - G(W_j^Y-)]^{n_1 \hat{p}_j^Y} \right)$$

$$= \left( \prod_{i=1}^{m_0} [F(W_i^X) - F(W_i^X-)]^{n_0 \hat{p}_i^X} \right)$$

$$\times \left( \prod_{j=1}^{m_1} \{\varphi(W_j^Y; \theta_0)[F(W_j^Y) - F(W_j^Y-)]\}^{n_1 \hat{p}_j^Y} \right).$$

Thus, from (3.2) the *weighted empirical likelihood function* for model (1.1) is given by

$$(3.3) \quad L(\theta, F) = \left( \prod_{i=1}^{m} p_i^{n\omega_i} \right) \left( \prod_{j=m_0+1}^{m} [\varphi(W_j; \theta)]^{n\omega_j} \right)$$

$$\text{for } p_i = F(W_i) - F(W_i-),$$

and the SPMLE $(\tilde{\theta}_n, \tilde{F}_n)$ for $(\theta_0, F_0)$ is the solution that maximizes $L(\theta, F)$. One may note that the use of weighted empirical likelihood function (2.1) here provides a simple and direct way to incorporate the model assumption of (1.1) in the derivation of likelihood function (3.3) for censored data. Also note that using the usual likelihood functions for specific types of censored data would result in a much more complicated likelihood function which is very difficult to handle.

To find $(\tilde{\theta}_n, \tilde{F}_n)$, we need to solve the following optimization problem:

$$\max L(\theta, \boldsymbol{p}) = \left( \prod_{i=1}^{m} p_i^{n\omega_i} \right) \left( \prod_{j=m_0+1}^{m} [\varphi(W_j; \theta)]^{n\omega_j} \right)$$

(3.4)
$$\text{subject to} \quad p_i \geq 0, \sum_{i=1}^{m} p_i = 1, \sum_{i=1}^{m} p_i \varphi(W_i; \theta) = 1,$$

where the last constraint reflects the fact that $\varphi(x;\theta)[F(x) - F(x-)]$ is a distribution function. Note that the NPMLE for censored data (1.2)–(1.7) is not always a proper d.f. (Mykland and Ren [17]). But for the moment, we assume $\sum_{i=1}^{m_0} \hat{p}_i^X = \sum_{i=1}^{m_1} \hat{p}_i^Y = 1$ in (3.1), which will not be needed later on for our main results of the paper. To solve (3.4), we first maximize $L(\theta, \boldsymbol{p})$ with respect to $\boldsymbol{p} = (p_1, \ldots, p_m)$ for fixed $\theta$, then maximize $l(\theta) = \ln L(\theta, \tilde{\boldsymbol{p}}) = \max_{\boldsymbol{p}} \ln L(\theta, \boldsymbol{p})$ over $\theta$ to find $\tilde{\theta}_n$. Noting that for $U_i(\theta) = \varphi(W_i; \theta)$, constraints in (3.4) imply $\sum_{i=1}^{m} p_i[U_i(\theta) - 1] = 0$, we know that $\theta$ must satisfy

(3.5) $$[U_{(1)}(\theta) - 1] < 0 < [U_{(m)}(\theta) - 1].$$

Using the Lagrange multiplier method, it can be shown that for any fixed $\theta$ satisfying (3.5), the convexity of $\ln L(\theta, \boldsymbol{p})$ ensures that $L(\theta, \boldsymbol{p})$ is uniquely maximized by $L(\theta, \tilde{\boldsymbol{p}})$ (see pages 90–91 and 164 of Bazaraa, Sherali and Shetty [1]), where

(3.6) $$\tilde{p}_i = \frac{\omega_i}{1 + \lambda(\theta)[U_i(\theta) - 1]}, \qquad i = 1, \ldots, m,$$

with $\lambda(\theta)$ as the unique solution on interval $(-[U_{(m)}(\theta) - 1]^{-1}, -[U_{(1)}(\theta) - 1]^{-1})$ for

(3.7) $$0 = \psi(\lambda; \theta) \equiv \sum_{i=1}^{m} \frac{\omega_i[U_i(\theta) - 1]}{1 + \lambda[U_i(\theta) - 1]}.$$

Thus, we have $l(\theta) = n \sum_{i=1}^{m} \omega_i \ln \tilde{p}_i + n \sum_{j=m_0+1}^{m} \omega_j \ln \varphi(W_j; \theta)$.

For our examples, we have $\theta_0 \in \mathbb{R}$ or $\theta_0 \in \mathbb{R}^2$ in (1.1), and that for some functions $h_1(\theta)$ and $h_2(x)$, the following assumption holds for $\varphi(x;\theta)$ with $\theta \in \mathbb{R}$ or $\theta \in \mathbb{R}^2$:

(AS0) $\nabla \varphi(x;\theta) = \varphi(x;\theta) h_1(\theta)(1, h_2(x))^\top$ for $\nabla = (\partial/\partial\theta_1, \partial/\partial\theta_2)^\top$,

where $0 < h_1(\theta) \in \mathbb{R}$ is twice differentiable for $\theta \in \Theta$; $0 \leq h_2(x) \in \mathbb{R}$ is monotone for $x \geq 0$; in the case $\theta \in \mathbb{R}$, we have degenerating $h_2(x) \equiv 0$; in the case $\theta \in \mathbb{R}^2$, we always have strictly monotone $h_2(x)$ on the support of $F_0$. Throughout this paper, our notations mean that for the case $\theta \in \mathbb{R}$, only the nondegenerating component in equations, vectors and matrices is meaningful. To minimize $l(\theta)$, from (3.2), (3.6)–(3.7), $\psi(\lambda(\theta); \theta) = 0$ and constraints in (3.4), we obtain that under assumption (AS0):

$$\frac{\partial l}{\partial \theta_1} = -n\lambda(\theta) h_1(\theta) \sum_{i=1}^{m} \tilde{p}_i \varphi(W_i; \theta) + n h_1(\theta) \sum_{j=m_0+1}^{m} \omega_j$$

$$(3.8) \qquad = nh_1(\theta)[\rho_1 - \lambda(\theta)],$$

$$\frac{\partial l}{\partial \theta_2} = nh_1(\theta)\left(\rho_1 \sum_{j=m_0+1}^{m} \hat{p}_j h_2(W_j) - \lambda(\theta) \sum_{i=1}^{m} \tilde{p}_i \varphi(W_i; \theta) h_2(W_i)\right),$$

where the use of $\nabla \lambda(\theta)$ in deriving (3.8) can easily be justified by the theorems on implicit functions in mathematical analysis. If $\tilde{\theta}_n$ is a solution of $\nabla l(\theta) = \mathbf{0}$, then

$$(3.9) \quad \lambda(\tilde{\theta}_n) = \rho_1 \quad \text{and} \quad \sum_{j=m_0+1}^{m} \hat{p}_j h_2(W_j) - \sum_{i=1}^{m} \tilde{p}_i \varphi(W_i; \tilde{\theta}_n) h_2(W_i) = 0.$$

In the Appendix, we show that $\tilde{\theta}_n$ is equivalently given by the solution of equation(s):

$$(3.10) \quad \begin{cases} 0 = g_1(\theta) \equiv \int_0^\infty \frac{\varphi(x;\theta)}{\rho_0 + \rho_1 \varphi(x;\theta)} d\hat{F}(x) - \int_0^\infty \frac{1}{\rho_0 + \rho_1 \varphi(x;\theta)} d\hat{G}(x), \\ 0 = g_2(\theta) \equiv \int_0^\infty \frac{\varphi(x;\theta)h_2(x)}{\rho_0 + \rho_1 \varphi(x;\theta)} d\hat{F}(x) - \int_0^\infty \frac{h_2(x)}{\rho_0 + \rho_1 \varphi(x;\theta)} d\hat{G}(x), \end{cases}$$

by which we always mean that $\tilde{\theta}_n \in \mathbb{R}$ is the solution of $g_1(\theta) = 0$ if $h_2(x) \equiv 0$. For our examples, the unique existence of solution $\tilde{\theta}_n$ for (3.10) is shown in Sections 4 and 5, respectively, and it can be shown that $\tilde{\theta}_n$ maximizes $l(\theta)$ over those $\theta$ satisfying (3.5) (the proofs are omitted). Thus, $\tilde{\theta}_n$ is the SPMLE for $\theta_0$ in (1.1). Consequently, replacing $\theta$ by $\tilde{\theta}_n$ in (3.6), we obtain the following SPMLE $\tilde{F}_n$ for $F_0$:

$$(3.11) \quad \tilde{F}_n(t) = \sum_{i=1}^{m} \tilde{p}_i I\{W_i \leq t\} = \int_0^t \frac{1}{\rho_0 + \rho_1 \varphi(x;\tilde{\theta}_n)} d[\rho_0 \hat{F}(x) + \rho_1 \hat{G}(x)].$$

Since the equations in (3.10) only depend on the NPMLE $\hat{F}$ and $\hat{G}$, thus for the rest of the paper, $\tilde{\theta}_n$ denotes the solution of (3.10) without assumption $\sum_{i=1}^{m_0} \hat{p}_i^X = \sum_{i=1}^{m_1} \hat{p}_i^Y = 1$ in (3.1), and is used to compute $\tilde{F}_n$ in (3.11). In the following theorems, some asymptotic results on $(\tilde{\theta}_n, \tilde{F}_n)$ are established under some of the assumptions listed below, while the proofs are deferred to the Appendix.

(AS1) (a) $\varphi(x;\theta)$ is monotone in $x$ for any fixed $\theta \in \Theta$, where $\Theta = \{\theta_1 | a_1 < \theta_1 < \infty\}$ if $\theta \in \mathbb{R}$; $\Theta = \{(\theta_1, \theta_2) | a_i < \theta_i < \infty, i = 1, 2\}$ if $\theta \in \mathbb{R}^2$;
  (b) $\varphi(x;\theta)$ is increasing in $\theta_1$ (and in $\theta_2$ if $\theta \in \mathbb{R}^2$) for any fixed $x > 0$;
  (c) for fixed $x > 0$ (and fixed $\theta_2$ if $\theta \in \mathbb{R}^2$), $\varphi(x;\theta) \to \infty(0)$, as $\theta_1 \to \infty(a_1)$;
  (d) for $\theta = (\theta_1, \theta_2) \in \mathbb{R}^2$ and fixed $x > 0$, when $-\theta_1/\theta_2 \to \gamma$ with $0 \leq \gamma \leq \infty$: $\varphi(x;\theta) \to 0(\infty)$ if $x < \gamma(x > \gamma)$, as $\theta_2 \to \infty$; $\varphi(x;\theta) \to 0(\infty)$ if $x > \gamma(x < \gamma)$, as $\theta_2 \to a_2$;

(AS2) $\rho_0 = \frac{n_0}{n}$ and $\rho_1 = \frac{n_1}{n}$ remain the same as $n \to \infty$;

(AS3) $\sqrt{n}_0 \int_0^\infty \frac{[h_2(x)]^k \varphi(x;\theta_0)}{\rho_0+\rho_1\varphi(x;\theta_0)} d[\hat{F}(x) - F_0(x)] \stackrel{D}{\to} N(0, \sigma_{F,k}^2)$, as $n \to \infty$, $\sqrt{n}_1 \times \int_0^\infty \frac{[h_2(x)]^k}{\rho_0+\rho_1\varphi(x;\theta_0)} d[\hat{G}(x) - G_0(x)] \stackrel{D}{\to} N(0, \sigma_{G,k}^2)$, as $n \to \infty$, where $k = 0, 1$, and $[h_2(x)]^0 \equiv 1$;

(AS4) $\|\hat{F} - F_0\| \stackrel{a.s.}{\to} 0$, $\|\hat{G} - G_0\| \stackrel{a.s.}{\to} 0$, as $n \to \infty$;

(AS5) $\int_0^\infty [h_2(x)]^k d[\hat{F}(x) - F_0(x)] \stackrel{a.s.}{\to} 0$, $\int_0^\infty [h_2(x)]^k d[\hat{G}(x) - G_0(x)] \stackrel{a.s.}{\to} 0$, as $n \to \infty$, with finite $\int_0^\infty [h_2(x)]^k dF_0(x)$ and $\int_0^\infty [h_2(x)]^k dG_0(x)$, where $k = 1, 2, 3$;

(AS6) $\sqrt{n}_0(\hat{F} - F_0) \stackrel{w}{\Rightarrow} \mathbb{G}_F$, $\sqrt{n}_1(\hat{G} - G_0) \stackrel{w}{\Rightarrow} \mathbb{G}_G$, as $n \to \infty$, where $\mathbb{G}_F$ and $\mathbb{G}_G$ are centered Gaussian processes.

THEOREM 1. *Assume* (AS0)–(AS5). *Under model (1.1), we have:*

(i) $\tilde{\theta}_n \stackrel{a.s.}{\to} \theta_0$, *as* $n \to \infty$;
(ii) $\sqrt{n}(\tilde{\theta}_n - \theta_0) \stackrel{D}{\to} N(0, \Sigma_0)$, *as* $n \to \infty$;
(iii) $\|\tilde{F}_n - F_0\| \stackrel{a.s.}{\to} 0$, *as* $n \to \infty$.

THEOREM 2. *Assume* (AS0)–(AS6). *Under model (1.1), we have that* $\sqrt{n}(\tilde{F}_n - F_0)$ *weakly converges to a centered Gaussian process.*

REMARK 1 (*Assumptions of theorems*). For our examples, (AS0)–(AS1) hold, which will be discussed in Sections 4 and 5, respectively. From Gill [7], Gu and Zhang [11], Huang [12], Huang and Wellner [13] and Geskus and Groeneboom [6], we know that under some suitable conditions, (AS3) holds for censored data (1.2)–(1.7) aforementioned. We also know that for these types of censored data, (AS4) holds under some suitable conditions; see Stute and Wang [30], Gu and Zhang [11], Huang [12] and Groeneboom and Wellner [10]. For right censored data, (AS5) holds under some regularity conditions (Stute and Wang [30]). For other types of censored data, (AS5) is implied by (AS4) if the support of $F_0$ is finite. On the other hand, if weaker consistency result is desired in Theorem 1(i), assumption (AS5) can be weakened. Moreover, from Gill [7], Gu and Zhang [11] and Huang [12], we know that (AS6) holds under some suitable conditions for right censored data, doubly censored data and partly interval-censored data. The techniques used in our proofs show that the weak convergence of $\tilde{F}_n$ for interval censored data relies on that of NPMLE $\hat{F}$ or $\hat{G}$ for interval censored data, which is now unknown.

**4. Biased sampling models.** For the biased sampling problem in Example 1, this section discusses assumptions (AS0)–(AS1), shows the unique existence of SPMLE $\tilde{\theta}_n$ for $\theta_0 \in \mathbb{R}$ in (1.8), and studies the weighted empirical log-likelihood ratio for $w_0$.

Under (1.8), we have that in (AS0), $h_1(\theta) = 1/\theta$ for $\theta \in \Theta = \{\theta | a_1 = 0 < \theta < \infty\}$ and $h_2(x) \equiv 0$, and that (AS1)(a)–(c) obviously hold for any monotone weight function $w(x)$, while (AS1)(d) does not apply. Since $h_2(x) \equiv 0$, $\tilde{\theta}_n \in \mathbb{R}$ is determined by the first equation of (3.10). Note that (AS1)(c) and the Dominated Convergence Theorem (DCT) imply: $\lim_{\theta \to 0} g_1(\theta) = -\hat{G}(\infty)/\rho_0 < 0$ and $\lim_{\theta \to \infty} g_1(\theta) = \hat{F}(\infty)/\rho_1 > 0$. Thus, the solution $\tilde{\theta}_n$ of equation $g_1(\theta) = 0$ uniquely exists because $g_1'(\theta) > 0$ for $\theta > 0$.

*Weighted empirical log-likelihood ratio.* From (3.3) and (3.6), we know that under (1.8), the weighted empirical likelihood ratio is given by $\hat{R}(F) = L(\theta, F)/L(\tilde{\theta}_n, \tilde{F}_n) = (\theta/\tilde{\theta}_n)^{n\rho_1} \prod_{i=1}^{m}(p_i/\tilde{p}_i)^{n\omega_i}$, where $F(x) = \sum_{i=1}^{m} p_i I\{W_i \leq x\}$, $\theta = 1/[\sum_{i=1}^{m} p_i w(W_i)]$ and $\tilde{p}_i = \omega_i/[\rho_0 + \rho_1 \tilde{\theta}_n w(W_i)]$. Then, set $S = \{\int w(x) \, dF(x) | \hat{R}(F) \geq c\}$ may be used as confidence interval for $w_0$, where $0 < c < 1$ is a constant. Let

$$(4.1) \quad r(\theta_0) = \sup\left\{ (\theta_0/\tilde{\theta}_n)^{n\rho_1} \prod_{i=1}^{m} (p_i/\tilde{p}_i)^{n\omega_i} \Big| p_i \geq 0, \sum_{i=1}^{m} p_i = 1, \sum_{i=1}^{m} p_i w(W_i) = \frac{1}{\theta_0} \right\}.$$

It is easy to show that $S$ is an interval expressed by $S = [X_L, X_U]$, and that $X_L \leq w_0 \leq X_U$ if and only if $r(\theta_0) \geq c$, where $X_L = \inf\{\int_0^\infty w(x) \, dF(x) | F \in \mathcal{F}\}$ and $X_U = \sup\{\int_0^\infty w(x) \, dF(x) | F \in \mathcal{F}\}$ for $\mathcal{F} = \{F | \hat{R}(F) \geq c, p_i \geq 0, \sum_{i=1}^{m} p_i = 1\}$. We call $[X_L, X_U]$ the *weighted empirical likelihood ratio confidence interval* for $w_0$, and the limiting distribution of weighted empirical log-likelihood ratio for those censored data (1.2)–(1.7) is given in the following theorem with a proof sketched in the Appendix.

THEOREM 3. *Assume* (AS2)–(AS5) *for model* (1.8). *Then,* $-2 \ln r(\theta_0) \xrightarrow{D} c_0 \chi_1^2$, *as* $n \to \infty$, *where* $0 < c_0 < \infty$ *is a constant and* $\chi_1^2$ *has a chi-squared distribution.*

**5. Case-control logistic regression models.** For the case-control logistic regression model in Example 2, this section discusses assumptions (AS0)–(AS1), shows the unique existence of SPMLE $\tilde{\theta}_n$ for $\theta_0 \in \mathbb{R}^2$ in (1.9), and provides a goodness-of-fit test for model (1.10).

Under (1.9), we have that in (AS0)–(AS1), $h_1(\theta) \equiv 1$ for $\theta \in \Theta$ with $a_1 = a_2 = -\infty$ and $h_2(x) = x$, and that (AS1) holds for $\varphi(x; \theta) = \exp(\alpha + \beta x)$ with $\theta = (\alpha, \beta) \in \mathbb{R}^2$. In the Appendix, we show that the solution $\tilde{\theta}_n$ of (3.10) exists uniquely.

*Goodness-of-fit test.* To assess the validity of logistic regression model assumption (1.10) with censored data, note that there are two ways to estimate d.f. $F_0$ in (1.9) using censored data (1.11). One is the NPMLE $\hat{F}$ based on the first sample, and the other is the SPMLE $\tilde{F}_n$ based on both samples under model assumption (1.10), that is, (1.9). Based on Theorems 1 and 2, we have the following corollary on the asymptotic properties of $\hat{F}$ and $\tilde{F}_n$ with proofs deferred to the Appendix.

COROLLARY 1. *Assume* (AS2)–(AS5) *for model (1.9). Then, as $n \to \infty$:*

(i) $\|\tilde{F}_n - \hat{F}\| \stackrel{\text{a.s.}}{\to} 0$ *under model (1.10);*

(ii) $\|\tilde{F}_n - F_1\| \stackrel{\text{a.s.}}{\to} 0$ *when model (1.10) does not hold [i.e., $g_0(x) \stackrel{\text{a.e.}}{=} \varphi(x; \theta_0) \times f_0(x)$ does not hold], where $F_1 \neq F_0$;*

(iii) $\sqrt{n}(\tilde{F}_n - \hat{F})$ *weakly converges to a centered Gaussian process under model (1.10) and assumption* (AS6).

Thus, from Remark 1 we know that for right censored data, doubly censored data and partly interval-censored data, we may use the following Kolmogorov–Smirnov-type statistic to measure the difference between $\hat{F}$ and $\tilde{F}_n$, which gives a goodness-of-fit test statistic for case-control logistic regression model (1.10):

$$\text{(5.1)} \qquad T_n = \sqrt{n}\|\tilde{F}_n - \hat{F}\| = \sqrt{n} \sup_{0 \leq t < \infty} |\tilde{F}_n(t) - \hat{F}(t)|.$$

*Bootstrap method.* To compute the $p$-value for test statistic $T_n$ in (5.1), we suggest the following $n$ out of $n$ bootstrap method. Since $\tilde{\theta}_n = (\tilde{\alpha}_n, \tilde{\beta}_n)$ is determined by (3.10), it is a functional of the NPMLE $\hat{F}$ and $\hat{G}$, denoted as $\tilde{\theta}_n = \theta(\hat{F}, \hat{G})$; in turn, (3.11) implies that $\tilde{F}_n(t) - \hat{F}(t)$ is a functional of $\hat{F}$ and $\hat{G}$, denoted as $\tilde{F}_n - \hat{F} = \tau(\hat{F}, \hat{G})$. Note that under model (1.1), $\theta_0$ is the unique solution of equation(s):

$$\text{(5.2)} \quad \begin{aligned} 0 &= g_{01}(\theta) \equiv \int_0^\infty \frac{\varphi(x;\theta)}{\rho_0 + \rho_1 \varphi(x;\theta)} dF_0(x) - \int_0^\infty \frac{1}{\rho_0 + \rho_1 \varphi(x;\theta)} dG_0(x), \\ 0 &= g_{02}(\theta) \equiv \int_0^\infty \frac{\varphi(x;\theta) h_2(x)}{\rho_0 + \rho_1 \varphi(x;\theta)} dF_0(x) - \int_0^\infty \frac{h_2(x)}{\rho_0 + \rho_1 \varphi(x;\theta)} dG_0(x), \end{aligned}$$

by which we always mean that $\theta_0 \in \mathbb{R}$ is the solution of $g_{01}(\theta) = 0$ if $h_2(x) \equiv 0$. Thus, under (1.9) we have $\theta_0 = (\alpha_0, \beta_0) = \theta(F_0, G_0)$; in turn, $\tau(F_0, G_0) \equiv 0$, which means $T_n = \sqrt{n}\|\tau(\hat{F}, \hat{G}) - \tau(F_0, G_0)\|$ under model (1.10). Hence, from the formulation given in Bickel and Ren [3], the distribution of $T_n$ under model (1.10) can be estimated by that of $T_n^* = \sqrt{n}\|\tau(\hat{F}^*, \hat{G}^*) - \tau(\hat{F}, \hat{G})\|$, where $\hat{F}^*$ and $\hat{G}^*$ are calculated based on the $n$ out of $n$ bootstrap samples,

respectively. For instance, $\hat{F}^*$ is calculated based on the bootstrap sample $\boldsymbol{O}_1^{X*}, \ldots, \boldsymbol{O}_{n_0}^{X*}$ taken with replacement from $\{\boldsymbol{O}_1^X, \ldots, \boldsymbol{O}_{n_0}^X\}$. The $p$-value is estimated by the percentage of $T_n^*$'s that are greater than test statistic $T_n$. Note that the $n$ out of $n$ bootstrap consistency for $\sqrt{n}_0(\hat{F} - F_0)$ estimated by $\sqrt{n}_0(\hat{F}^* - \hat{F})$ has been established for right censored data, doubly censored data and partly interval-censored data by Bickel and Ren [2] and Huang [12].

REMARK 2. The proposed test (5.1) can be used for any type of censored data as long as (AS2)–(AS6) hold. When (AS6) does not hold, such as for interval censored data, Corollary 1 shows that we may graphically check the model fitting for (1.10) by comparing curves of $\hat{F}$ and $\tilde{F}_n$. Note that when model (1.10) does not hold, statistic $T_n^*$ is still asymptotically a function of a centered Gaussian process, but $T_n \overset{\text{a.s.}}{\to} \infty$ based on Corollary 1(ii). Thus, our proposed test is consistent. In terms of computing $(\tilde{\alpha}_n, \tilde{\beta}_n)$, it can be done using the Newton–Raphson method described on page 374 of Press et al. [21] to solve (3.10); a computation routine in FORTRAN is available from the author. Although not presented here, our extensive simulation studies on $(\tilde{\alpha}_n, \tilde{\beta}_n)$ and the comparison between the distributions of $T_n$ and $T_n^*$ give excellent results.

## APPENDIX

PROOF OF "$\tilde{\theta}_n$ IS EQUIVALENTLY GIVEN BY THE SOLUTION OF (3.10)." Under assumption $\sum_{i=1}^{m_0} \hat{p}_i^X = \sum_{i=1}^{m_1} \hat{p}_i^Y = 1$, the first equation of (3.9) is equivalent to $\psi(\rho_1; \theta) = 0$, which by (3.7) and (3.1)–(3.2), gives $g_1(\theta) = 0$ in (3.10). The proof follows from that (3.6) and $\lambda(\theta) = \rho_1$ imply that the second equation of (3.9) is $0 = -\rho_0 g_2(\theta)$. □

PROOF OF "UNIQUE EXISTENCE OF $\tilde{\theta}_n$ IN EXAMPLE 2." Let

$$(A.1) \quad R_n(\theta) = \int_0^\infty \frac{1}{\rho_1} \ln[\rho_0 + \rho_1 \varphi(x; \theta)] \, d\hat{F}(x) + \int_0^\infty \frac{1}{\rho_0} \ln\left(\frac{\rho_0 + \rho_1 \varphi(x; \theta)}{\varphi(x; \theta)}\right) d\hat{G}(x).$$

Since $\hat{F}$ and $\hat{G}$ are step functions with finite jumps, we know that $R_n(\theta)$ is well defined on $\mathbb{R}^2$. From (A.1) and (3.10), we have $\nabla R_n(\theta) = h_1(\theta)(g_1(\theta), g_2(\theta))^\top$ and

$$\boldsymbol{\Sigma}_{R_n, \theta} = \begin{pmatrix} \dfrac{\partial^2 R_n}{\partial \theta_1^2} & \dfrac{\partial^2 R_n}{\partial \theta_2 \, \partial \theta_1} \\ \dfrac{\partial^2 R_n}{\partial \theta_1 \, \partial \theta_2} & \dfrac{\partial^2 R_n}{\partial \theta_2^2} \end{pmatrix}$$

$$
\begin{aligned}
\text{(A.2)} \quad &= (g_1(\theta), g_2(\theta))^\top (\nabla h_1(\theta))^\top \\
&\quad + h_1^2(\theta) \int_0^\infty \begin{pmatrix} 1 & h_2(x) \\ h_2(x) & h_2^2(x) \end{pmatrix} \\
&\qquad\qquad \times \frac{\varphi(x;\theta)}{[\rho_0 + \rho_1 \varphi(x;\theta)]^2} \, d[\rho_0 \hat{F}(x) + \rho_1 \hat{G}(x)].
\end{aligned}
$$

Thus, $\nabla R_n(\theta) = 0$ is equivalent to (3.10) because $h_1(\theta) > 0$ by (AS0). For Example 2, we have $h_1(\theta) \equiv 1$ and $h_2(x) = x$, which imply that $\Sigma_{R_n,\theta}$ is a positive-definite matrix. Hence, $R_n(\theta)$ is strictly convex. Moreover, note that under (1.9), we have in (A.1) $R_n(\theta) \geq (\ln \rho_0)/\rho_1 + (\ln \rho_1)/\rho_0$ for any $\theta = (\alpha, \beta) \in \mathbb{R}^2$, and that by a similar argument used in (6.5) of Ren and Gu [27], we can show: $\lim_{\lambda \to \infty} \inf R_n(\lambda e_1, \lambda e_2) = \infty$ for any $e_1^2 + e_2^2 = 1$. Hence, $R_n(\theta)$ has a unique global minimum point which must be the solution of (3.10) (see pages 101–102 of Bazaraa, Sherali and Shetty [1]).  □

PROOF OF THEOREM 1(i). Let $\hat{\mu}(x) = \rho_0 \hat{F}(x) + \rho_1 \hat{G}(x)$; then (3.10) gives

$$
\text{(A.3)} \quad
\begin{aligned}
\hat{F}(\infty) &= \int_0^\infty \frac{d\hat{\mu}(x)}{\rho_0 + \rho_1 \varphi(x;\tilde{\theta}_n)} \leq \frac{1}{\rho_0}, \\
\hat{G}(\infty) &= \int_0^\infty \frac{\varphi(x;\tilde{\theta}_n)\, d\hat{\mu}(x)}{\rho_0 + \rho_1 \varphi(x;\tilde{\theta}_n)} \leq \frac{1}{\rho_1},
\end{aligned}
$$

where (AS4) implies $\hat{F}(\infty) \overset{\text{a.s.}}{\to} 1, \hat{G}(\infty) \overset{\text{a.s.}}{\to} 1$, as $n \to \infty$. As follows, we show $\tilde{\theta}_n = O(1)$ almost surely for case $\tilde{\theta}_n = (\tilde{\theta}_n^{(1)}, \tilde{\theta}_n^{(2)}) \in \mathbb{R}^2$ (the proof for case $\tilde{\theta}_n \in \mathbb{R}$ is similar).

Assume $\tilde{\theta}_n^{(2)} \geq 0$. If $\tilde{\theta}_n^{(1)} \to \infty$, then from integration by parts, the boundedness of the integrand function, (AS1)(b)–(c) and the DCT, we have that in (A.3):

$$
\text{(A.4)} \quad 1 = \lim_{n \to \infty} \int_0^\infty \frac{d\mu_0(x)}{\rho_0 + \rho_1 \varphi(x;\tilde{\theta}_n)} \leq \int_0^\infty \lim_{n \to \infty} \frac{d\mu_0(x)}{\rho_0 + \rho_1 \varphi(x;\tilde{\theta}_n^{(1)}, 0)} = 0,
$$

a contradiction, where $\mu_0(x) = \rho_0 F_0(x) + \rho_1 G_0(x)$. Thus, $\tilde{\theta}_n^{(2)} \geq 0$ implies $\tilde{\theta}_n^{(1)} = O(1)$ or $\tilde{\theta}_n^{(1)} \to a_1$. Similarly, we know that $0 \leq \tilde{\theta}_n^{(2)} \leq M_2 < \infty$ and $\tilde{\theta}_n^{(1)} \to a_1$ imply $1 = \lim \int_0^\infty [\rho_0 + \rho_1 \varphi(x;\tilde{\theta}_n)]^{-1} d\mu_0(x) \geq \int_0^\infty \lim [\rho_0 + \rho_1 \varphi(x;\tilde{\theta}_n^{(1)}, M_2)]^{-1} d\mu_0(x) = 1/\rho_0$, a contradiction. Hence, if $\tilde{\theta}_n^{(2)} \geq 0$, then $\tilde{\theta}_n^{(2)} = O(1)$ implies $\tilde{\theta}_n^{(1)} = O(1)$.

Assume $\tilde{\theta}_n^{(2)} \to \infty, -\tilde{\theta}_n^{(1)}/\tilde{\theta}_n^{(2)} \to \gamma$ with $0 \leq \gamma \leq \infty$. Similarly as (A.4), (AS1) gives

$$
\text{(A.5)} \quad 1 = \int_0^\infty \lim_{n \to \infty} \frac{d\mu_0(x)}{\rho_0 + \rho_1 \varphi(x;\tilde{\theta}_n)} = \frac{\mu_0(\gamma)}{\rho_0},
$$

where we must have $0 < \gamma < \infty$ to be inside the support of $F_0$; a contradiction otherwise. Also, if we let $n \to \infty$ in the second equation of (3.10), then from (AS4)–(AS5), Hölder's inequality, the DCT and an argument similar to above, we have

$$\text{(A.6)} \qquad \frac{1}{\rho_1} \int_\gamma^\infty h_2(x)\, dF_0(x) = \frac{1}{\rho_0} \int_0^\gamma h_2(x)\, dG_0(x).$$

However, (A.5)–(A.6) contradict $[G_0(\gamma) \int_\gamma^\infty h_2(x)\, dF_0(x) - \bar{F}_0(\gamma) \times \int_0^\gamma h_2(x)\, dG_0(x)] = \iint_{x<\gamma<y}[h_2(y) - h_2(x)]\, dF_0(y)\, dG_0(x) \neq 0$, which is implied by (AS0). Thus, if $\tilde{\theta}_n^{(2)} \geq 0$, we must have $\tilde{\theta}_n^{(2)} = O(1)$; in turn, $\tilde{\theta}_n^{(1)} = O(1)$. Similarly, we can show $\tilde{\theta}_n^{(2)} = O(1)$ and $\tilde{\theta}_n^{(1)} = O(1)$ if $\tilde{\theta}_n^{(2)} < 0$. Hence, we have $\tilde{\theta}_n = O(1)$ almost surely.

Assume $\tilde{\theta}_n \to \eta_0$, as $n \to \infty$. Then, from (3.10) and an argument similar to that used in (A.6), we know that $\eta_0$ is a solution of (5.2). Note that for nondegenerating $h_2(x)$, to obtain the second equation of (5.2) for $\eta_0$ we use (AS5) and the proof of Lemma 3 of Gill [8], noticing that $h_2(x)$ is monotone and $[\rho_0 + \rho_1 \varphi(x; \eta_0)]^{-1}$ is bounded and continuous. Hence, the proof follows from the uniqueness of the solution for (5.2). $\square$

PROOF OF THEOREM 1(ii). Here, we only prove the case $\tilde{\theta}_n \in \mathbb{R}^2$, because the proof for case $\tilde{\theta}_n \in \mathbb{R}$ is similar. For $R_n(\theta)$ in (A.1), we have that under model (1.1):

$$\text{(A.7)} \qquad \begin{aligned} \nabla R_n(\theta_0) &= h_1(\theta_0)([g_1(\theta_0) - g_{01}(\theta_0)], [g_2(\theta_0) - g_{02}(\theta_0)])^\top, \\ \nabla R_n(\tilde{\theta}_n) &= \nabla R_n(\theta_0) + \boldsymbol{\Sigma}_{R_n,\theta_0}(\tilde{\theta}_n - \theta_0)^\top + \tfrac{1}{2}(r_1(\tilde{\theta}_n), r_2(\tilde{\theta}_n))^\top, \end{aligned}$$

where $g_1, g_2$ and $g_{01}, g_{02}$ are given in (3.10) and (5.2), respectively; $\boldsymbol{\Sigma}_{R_n,\theta}$ is given in (A.2); and from (AS5), Theorem 1(i) and straightforward calculation based on (A.2), we have $r_i(\tilde{\theta}_n) = o_p(\tilde{\theta}_n - \theta_0)$. From (A.7), (AS3), the independence between $\hat{F}$ and $\hat{G}$, and page 4 of Serfling [29], we know that $\sqrt{n} \nabla R_n(\theta_0)$ converges in distribution to a normal random vector, while (A.2), (5.2) and a similar argument in (A.6) imply

$$\text{(A.8)} \qquad \boldsymbol{\Sigma}_{R_n,\theta_0} \stackrel{\text{a.s.}}{\to} \boldsymbol{\Sigma}_1 = h_1^2(\theta_0) \int_0^\infty \begin{pmatrix} 1 & h_2(x) \\ h_2(x) & h_2^2(x) \end{pmatrix} \frac{\varphi(x; \theta_0)\, d\mu_0(x)}{[\rho_0 + \rho_1 \varphi(x; \theta_0)]^2}$$

as $n \to \infty$,

where $\boldsymbol{\Sigma}_1$ is positive-definite. Hence, $\nabla R_n(\tilde{\theta}_n) = 0$, (A.7)–(A.8) and Theorem 1(i) give

$$\text{(A.9)} \qquad \sqrt{n}(\tilde{\theta}_n - \theta_0) = -\boldsymbol{\Sigma}_1^{-1} \sqrt{n} \nabla R_n(\theta_0) + o_p(1). \qquad \square$$

PROOF OF THEOREM 1(iii). Here, we only prove the case $\tilde{\theta}_n \in \mathbb{R}^2$, because the proof for case $\tilde{\theta}_n \in \mathbb{R}$ is similar. For any $t > 0$, we let $\tilde{F}_n(t) \equiv g_3(\tilde{\theta}_n)$ in (3.11); then

$$\text{(A.10)} \quad \begin{aligned} \tilde{F}_n(t) &= g_3(\tilde{\theta}_n) \\ &= g_3(\theta_0) + (\tilde{\theta}_n - \theta_0)\nabla g_3(\theta_0) + \tfrac{1}{2}(\tilde{\theta}_n - \theta_0)\boldsymbol{\Sigma}_{g_3,\xi_n}(\tilde{\theta}_n - \theta_0)^\top, \end{aligned}$$

where $\xi_n$ is between $\tilde{\theta}_n$ and $\theta_0$, and

$$\nabla g_3(\theta) = -\rho_1 h_1(\theta) \int_0^t (1, h_2(x))^\top \frac{\varphi(x;\theta)}{[\rho_0 + \rho_1 \varphi(x;\theta)]^2} \, d\hat{\mu}(x),$$

$$\boldsymbol{\Sigma}_{g_3,\theta} = \begin{pmatrix} \frac{\partial^2 g_3}{\partial \theta_1^2} & \frac{\partial^2 g_3}{\partial \theta_2 \partial \theta_1} \\ \frac{\partial^2 g_3}{\partial \theta_1 \partial \theta_2} & \frac{\partial^2 g_3}{\partial \theta_2^2} \end{pmatrix}$$

$$\text{(A.11)} \quad \begin{aligned} &= [h_1(\theta)]^{-1} \nabla g_3(\theta) [\nabla h_1(\theta)]^\top \\ &\quad - \rho_1 h_1^2(\theta) \int_0^t \begin{pmatrix} 1 & h_2(x) \\ h_2(x) & h_2^2(x) \end{pmatrix} \\ &\quad \times \frac{\varphi(x;\theta)[\rho_0 - \rho_1 \varphi(x;\theta)]}{[\rho_0 + \rho_1 \varphi(x;\theta)]^3} \, d\hat{\mu}(x). \end{aligned}$$

From (AS5) and Theorem 1(i), we have that uniformly in $t$,

$$\left| \frac{\partial^2 g_3(\xi_n)}{\partial \theta_2^2} \right| \leq \rho_1 \int_0^\infty \frac{\varphi(x;\xi_n) h_2(x)[(\partial^2 h_1(\xi_n)/\partial \theta_2^2) + h_1^2(\xi_n)h_2(x)]}{[\rho_0 + \rho_1 \varphi(x;\xi_n)]^2} \, d\hat{\mu}(x)$$
$$= O_{\text{a.s.}}(1),$$

which also holds for other partial derivatives in (A.11). Thus, Theorem 1(ii) implies that with $(\tilde{\theta}_n - \theta_0)\nabla g_3(\theta_0) = o_{\text{a.s.}}(1)$, (A.10) can be written as

$$\text{(A.12)} \quad \tilde{F}_n(t) = g_3(\theta_0) + (\tilde{\theta}_n - \theta_0)\nabla g_3(\theta_0) + O_{\text{a.s.}}(|\tilde{\theta}_n - \theta_0|^2).$$

From (AS4) and integration by parts, we have $|g_3(\theta_0) - F_0(t)| \overset{\text{a.s.}}{\to} 0$ for any fixed $t > 0$; in turn, the proof follows from (A.12) and Pólya's Theorem. □

PROOF OF THEOREM 2. Here, we only prove the case $\tilde{\theta}_n \in \mathbb{R}^2$, because the proof for case $\tilde{\theta}_n \in \mathbb{R}$ is similar. Let $(\hat{v}_1, \hat{v}_2)^\top = \nabla g_3(\theta_0)$ as in (A.11), and let $(v_1, v_2)^\top = -\rho_1 h_1(\theta_0) \int_0^t (1, h_2(x))^\top \varphi(x;\theta_0)[\rho_0 + \rho_1 \varphi(x;\theta_0)]^{-2} \, d\mu_0(x)$. From (AS4) and integration by parts, we have $|\hat{v}_k(t) - v_k(t)| \overset{\text{a.s.}}{\to} 0$ for any fixed $t > 0$, where $k = 1, 2$. Since $\hat{v}_k(t)$ and $v_k(t)$ are continuous and monotone in $t$, then from (AS5) and a similar argument used in the proof of Theorem 1(i) for showing $\eta_0$ as the solution of (5.2), we have $\|\hat{v}_k - v_k\| \overset{\text{a.s.}}{\to} 0$, as $n \to \infty$. Thus, if we let $u_0(x) = 1/[\rho_0 + \rho_1 \varphi(x; \theta_0)]$, $u_1(x) = u_0(x)\varphi(x; \theta_0)$, and $\lambda_{ij}$ the elements of $\mathbf{\Sigma}_1^{-1}$, then (A.7), (A.9), (A.12), (1.1) and Theorem 1(ii) imply

$$\sqrt{n}[\tilde{F}_n(t) - F_0(t)]$$
$$= o_p(1)$$
(A.13)
$$+ \sqrt{n}\left(\int_0^t u_0(x)\,d\hat{\mu}(x) - F_0(t) - (\hat{v}_1(t), \hat{v}_2(t))\mathbf{\Sigma}_1^{-1}\nabla R_n(\theta_0)\right)$$
$$= o_p(1) + \sqrt{n}(\hat{U}_F - U_F) + \sqrt{n}(\hat{U}_G - U_G),$$

where for $s_1(t) = h_1(\theta_0)[\lambda_{11}v_1(t) + \lambda_{21}v_2(t)]$ and $s_2(t) = h_1(\theta_0)[\lambda_{12}v_1(t) + \lambda_{22}v_2(t)]$,

$$\sqrt{n}[\hat{U}_F(t) - U_F(t)]$$
$$\equiv \tau_1(\sqrt{n_0}(\hat{F} - F_0))$$
(A.14)
$$= \sqrt{n}\bigg(\rho_0 \int_0^t u_0(x)\,d[\hat{F}(x) - F_0(x)]$$
$$- s_1(t) \int_0^\infty u_1(x)\,d[\hat{F}(x) - F_0(x)]$$
$$- s_2(t) \int_0^\infty u_1(x)h_2(x)\,d[\hat{F}(x) - F_0(x)]\bigg),$$

$$\sqrt{n}[\hat{U}_G(t) - U_G(t)]$$
$$\equiv \tau_2(\sqrt{n_1}(\hat{G} - G_0))$$
(A.15)
$$= \sqrt{n}\bigg(\rho_1 \int_0^t u_0(x)\,d[\hat{G}(x) - G_0(x)]$$
$$+ s_1(t) \int_0^\infty u_0(x)\,d[\hat{G}(x) - G_0(x)]$$
$$+ s_2(t) \int_0^\infty u_0(x)h_2(x)\,d[\hat{G}(x) - G_0(x)]\bigg).$$

As (A.14) is a linear functional of $\sqrt{n_0}(\hat{F} - F_0)$, (AS6) implies $\sqrt{n}(\hat{U}_F - U_F) \stackrel{w}{\Rightarrow} \tau_1(\mathbb{G}_F)$, as $n \to \infty$, where from pages 154–157 of Iranpour and Chacon [14], we know that $\tau_1(\mathbb{G}_F)$ is a centered Gaussian process. Similarly, $\sqrt{n}[\hat{U}_G - U_G]$ in (A.15) weakly converges to a centered Gaussian process $\tau_2(\mathbb{G}_G)$. The proof follows from (A.13)–(A.15), and that $\tau_1(\mathbb{G}_F)$ and $\tau_2(\mathbb{G}_G)$ are two independent centered Gaussian processes. $\square$

PROOF OF COROLLARY 1. Note that part (i) follows directly from Theorem 1(iii) and (AS4), while part (iii) follows from some minor adjustments in the proof of Theorem 2. Thus, we only give the proof of part (ii) as follows.

Here, we have $h_1(\theta) \equiv 1$ and $h_2(x) = x$; thus in (A.2) we have $\nabla h_1(\theta) \equiv 0$. From the proofs of the unique existence of $\tilde{\theta}_n$ and Theorem 1(i), we know that when model (1.10) does not hold, $\tilde{\theta}_n$ is still well defined, and satisfies $|\tilde{\theta}_n - \theta_1| \stackrel{\text{a.s.}}{\to} 0$, as $n \to \infty$, where $\theta_1 = (\alpha_1, \beta_1)$ is the unique solution of (5.2) for $\varphi(x; \theta) = \exp(\alpha + \beta x)$. Applying this, (AS4) and integration by parts to (3.11), we have $\|\tilde{F}_n - F_1\| \stackrel{\text{a.s.}}{\to} 0$, where $F_1(t) = \int_0^t [\rho_0 + \rho_1 \exp(\alpha_1 + \beta_1 x)]^{-1} d\mu_0(x)$. It is easy to verify that $F_1 \neq F_0$ when (1.10) does not hold [otherwise, we have $g_0(x) \stackrel{\text{a.e.}}{=} \varphi(x; \theta_1) f_0(x)$ with $\theta_1 = \theta_0$], and that the first equation of (5.2) implies that $F_1$ is a distribution function. $\square$

PROOF OF THEOREM 3. For a simpler argument, we assume $\sum_{i=1}^{m_0} \hat{p}_i^X = \sum_{i=1}^{m_1} \hat{p}_i^Y = 1$ in (3.1), which can be removed with some additional work in our proof here. To get an expression of $r(\theta_0)$, it can be shown by using the Lagrange multiplier method that the solution of the maximization problem in (4.1) is $\bar{p}_i = \omega_i/(1 + \lambda_0 U_i)$, $1 \leq i \leq m$, where $U_i = [\theta_0 w(W_i) - 1]$, $\omega_i$ is given in (3.2), and $\lambda_0$ is the unique solution of equation $\phi(\lambda) = 0$ on interval $(-U_{(m)}^{-1}, -U_{(1)}^{-1})$ for $\phi(\lambda) \equiv \sum_{i=1}^{m} \bar{p}_i U_i = \sum_{i=1}^{m} (\omega_i U_i)/(1 + \lambda U_i)$. Thus, we have

$$(\text{A.16}) \quad \ln r(\theta_0) = -n \sum_{i=1}^{m} \omega_i \ln\left(\frac{1 + \lambda_0 U_i}{\rho_0 + \rho_1 \tilde{\theta}_n w(W_i)}\right) - n\rho_1 \ln(\tilde{\theta}_n/\theta_0).$$

Using Taylor's expansion on $\phi(\lambda)$, we have that from $\psi(\rho_1; \tilde{\theta}_n) = 0$ in (3.7),

$$(\text{A.17}) \quad \begin{aligned} \phi'(\xi)(\rho_1 - \lambda_0) &= \phi(\rho_1) - \phi(\lambda_0) \\ &= \sum_{i=1}^{m} (\omega_i U_i)/(1 + \rho_1 U_i) \\ &= \sum_{i=1}^{m} \frac{\omega_i[\theta_0 w(W_i) - 1]}{\rho_0 + \rho_1 \theta_0 w(W_i)} - \sum_{i=1}^{m} \frac{\omega_i[\tilde{\theta}_n w(W_i) - 1]}{\rho_0 + \rho_1 \tilde{\theta}_n w(W_i)} \\ &= \sum_{i=1}^{m} \frac{\omega_i w(W_i)(\theta_0 - \tilde{\theta}_n)}{[\rho_0 + \rho_1 \xi_i w(W_i)]^2}, \end{aligned}$$

where $\xi$ is between $\rho_1$ and $\lambda_0$, and $\xi_i$ is between $\theta_0$ and $\tilde{\theta}_n$. From (AS4), integration by parts and Theorem 1(i), we know that (A.17) implies

$$(\text{A.18}) \ (\rho_1 - \lambda_0) = (\theta_0 - \tilde{\theta}_n)\left(\frac{1}{\phi'(\rho_1)} \sum_{i=1}^{m} \frac{\omega_i w(W_i)}{[\rho_0 + \rho_1 \theta_0 w(W_i)]^2} + o_p(1)\right).$$

Also using Taylor's expansion, we have

$$\ln[\rho_0 + \rho_1 \theta_0 w(W_i)]$$

$$(\text{A.19}) \quad = \ln(1 + \rho_1 U_i) = \ln(1 + \lambda_0 U_i) + \frac{U_i}{1 + \lambda_0 U_i}(\rho_1 - \lambda_0)$$
$$- \frac{U_i^2}{2(1 + \lambda_0 U_i)^2}(\rho_1 - \lambda_0)^2 + \frac{U_i^3}{6(1 + \eta_i U_i)^3}(\rho_1 - \lambda_0)^3,$$

$$\ln[\rho_0 + \rho_1 \theta_0 w(W_i)]$$
$$(\text{A.20}) \quad = \ln[\rho_0 + \rho_1 \tilde{\theta}_n w(W_i)] + \frac{\rho_1 w(W_i)}{\rho_0 + \rho_1 \tilde{\theta}_n w(W_i)}(\theta_0 - \tilde{\theta}_n)$$
$$- \frac{[\rho_1 w(W_i)]^2}{2[\rho_0 + \rho_1 \tilde{\theta}_n w(W_i)]^2}(\theta_0 - \tilde{\theta}_n)^2$$
$$+ \frac{[\rho_1 w(W_i)]^3}{6[\rho_0 + \rho_1 \zeta_i w(W_i)]^3}(\theta_0 - \tilde{\theta}_n)^3,$$

$$(\text{A.21}) \quad -\ln(\tilde{\theta}_n/\theta_0)$$
$$= \frac{\theta_0 - \tilde{\theta}_n}{\tilde{\theta}_n} - \frac{(\theta_0 - \tilde{\theta}_n)^2}{2\tilde{\theta}_n^2} + \frac{(\theta_0 - \tilde{\theta}_n)^3}{6\zeta},$$

where $\eta_i$ is between $\rho_1$ and $\lambda_0$, while $\zeta_i$ and $\zeta$ are between $\theta_0$ and $\tilde{\theta}_n$. Since (A.18) and Theorem 1(ii) imply $(\rho_1 - \lambda_0) = O_p(n^{-1/2})$, then from $\sum_{i=1}^m (\omega_i U_i)/(1 + \lambda_0 U_i) = 0$ and $\tilde{\theta}_n \sum_{i=1}^m [\omega_i w(W_i)]/[\rho_0 + \rho_1 \tilde{\theta}_n w(W_i)] = \sum_{i=1}^m \tilde{p}_i \varphi(W_i; \tilde{\theta}_n) = 1$, and by applying (A.19)–(A.21) to (A.16), we obtain

$$(\text{A.22}) \quad \ln r(\theta_0) = O_p(n^{-1/2}) - \frac{n(\rho_1 - \lambda_0)^2}{2} \sum_{i=1}^m \frac{\omega_i U_i^2}{(1 + \lambda_0 U_i)^2}$$
$$- \frac{n(\tilde{\theta}_n - \theta_0)^2}{2}\left(\frac{\rho_1}{\tilde{\theta}_n^2} - \sum_{i=1}^m \frac{\omega_i[\rho_1 w(W_i)]^2}{[\rho_0 + \rho_1 \tilde{\theta}_n w(W_i)]^2}\right).$$

Hence, the proof follows from Theorem 1(ii) and applying (A.18) to (A.22), where the limits of the coefficients of $n(\tilde{\theta}_n - \theta_0)^2$ are handled similarly to (A.8). □

**Acknowledgments.** The author is very grateful to the referees, the Associate Editor and Editor Jianqing Fan for their comments and suggestions on the earlier version of the manuscript, which led to a much improved paper.

DEPARTMENT OF MATHEMATICS
UNIVERSITY OF CENTRAL FLORIDA
ORLANDO, FLORIDA 32816
USA
E-MAIL: jren@mail.ucf.edu